\documentstyle[12pt,amsmath,amsfonts,amssymb]{article}
\newtheorem{theorem}{Theorem}[section]

\newtheorem{proposition}[theorem]{Proposition}
\newtheorem{remark}[theorem]{Remark}
\newtheorem{lemma}[theorem]{Lemma}
\newcommand {\Mc}      {{\cal M}}

\newcommand {\Pc}      {{\cal P}}

\newcommand {\Tc}      {{\cal T}}
\newcommand {\Gc}      {{\cal G}}

\newcommand {\R}       {{\bf R}}
\newcommand {\RN}      {\R^n}

\newcommand {\rt}      {\rho_{\Tr}}
\newcommand {\dtr}     {{\rm deg}_{\Tr}}
\newcommand {\COR}     {C^{k,\omega}(\R)}
\newcommand {\CKH}     {\dot{C}^{k,\omega }(\R^{n})}
\newcommand {\CKO}     {C^{k,\omega }(\R^{n})}
\newcommand {\COM}     {C^{1,\omega }(\R^{n})}

\newcommand {\PK}      {\Pc_{k}}
\newcommand {\TK}      {\Tc_{k,n}}
\newcommand {\PRN}     {\PK\times\RN}
\newcommand {\SO}      {S_\omega}
\newcommand {\ro}      {r_\omega}
\newcommand {\MR}      {(\Mc,\rho)}
\newcommand {\card}    {\operatorname{card}}
\newcommand {\Lip}     {\operatorname{Lip}}
\newcommand {\BLip}    {\operatorname{{\bf Lip}}}

\newcommand {\Or}      {\operatorname{Or}}

\newcommand {\dd}      {\operatorname{d}}

\newcommand {\Tr}      {\operatorname{T}}
\newcommand {\emp}     {\emptyset}

\newcommand {\ts}      {\theta}
\newcommand {\vf}      {\varphi}
\newcommand {\vfa}     {\vf_{\alpha}}
\newcommand {\dm}      {\dd_\omega}
\newcommand {\dom}     {\delta_{\omega}}
\newcommand {\tP}      {\widetilde{P}}
\newcommand {\tT}      {\widetilde{T}}
\newcommand {\twG}     {\widetilde{\Gc}}
\newcommand {\tg}      {\tilde{g}}
\newcommand {\ty}      {\tilde{y}}

\newcommand {\ELG}      {\ell_{\Gc}}
\newcommand {\nn}      {\nonumber}
\newcommand {\rf}[1]   {(\ref{#1})}

\newcommand {\bx}      {$\hfill\Box$}
\newcommand {\SECT}[2] {\section*{\centerline{\normalsize
{\bf #1}}} \setcounter{section}{#2}
\setcounter{theorem}{0}\setcounter{equation}{0}}
\newcommand{\be}       {\begin{eqnarray}}
\newcommand{\bel}[1]   {\begin{eqnarray}\label{#1}}
\newcommand{\ee}       {\end{eqnarray}}
\begin{document}
\medskip
\centerline{\large{\bf The Whitney extension problem and }}
\vspace*{8mm}
\centerline{\large{\bf  Lipschitz selections of set-valued
mappings in jet-spaces}}
\vspace*{8mm}
\centerline{By {\it Pavel Shvartsman}} \vspace*{5 mm}
\centerline {\it Department of Mathematics, Technion -
Israel Institute of Technology,}\vspace*{2 mm}
\centerline{\it 32000 Haifa, Israel}\vspace*{2 mm}
\centerline{\it e-mail: pshv@tx.technion.ac.il}
\vspace*{10 mm}
\renewcommand{\thefootnote}{ }
\footnotetext[1]{{\it\hspace{-6mm}Math Subject
Classification} 46E35, 49K24, 52A35, 54C60, 54C65 \\
{\it Key Words and Phrases}~ Whitney's extension problem,
smooth functions, finiteness, metric, jet-space, set-valued
mapping, Lipschitz selection}
\begin{abstract}
{\small
\par We study a variant of the Whitney extension problem
\cite{W1,W2} for the space $\CKO$. We identify $\CKO$ with
a space of {\it Lipschitz} mappings from $\RN$ into the
space $\PRN$ of polynomial fields on $\RN$ equipped with a
certain metric. This identification allows us to
reformulate the Whitney problem for $\CKO$ as a Lipschitz
selection problem for set-valued mappings into a certain
family of subsets of $\PRN$. We prove a Helly-type
criterion for the existence of Lipschitz selections for
such set-valued mappings defined on finite sets. With the
help of this criterion, we improve estimates for finiteness
numbers in finiteness theorems for $\CKO$ due to C.
Fefferman \cite{F1,F3,F4}.
}
\end{abstract}
\vspace*{15mm}
\renewcommand{\thefootnote}{\arabic{footnote}}
\setcounter{footnote}{0}
\SECT{1. The main problem and main results}{1}
\indent
\par Let $\omega :\R_{+}\to \R_+$ be a continuous concave
function satisfying $\omega (0)=0$. We let $\CKH$ denote
the (homogeneous) space of all functions $f:\RN\to\R$ with
continuous derivatives of all orders up to $k$, for which
the seminorm
$$
\|f\|_{\CKH}:=\sum_{|\alpha|=k} \sup_{x,y\in\RN,x\neq y}
\frac{|D^{\alpha }f(x)-D^{\alpha}f(y)|}{\omega(\|x-y\|)}
$$
is finite. By $\CKO$ we denote the Banach subspace of
$\CKH$ defined by the norm
$$
\|f\|_{\CKO}:=\sum_{|\alpha|\le k} \sup_{x\in\RN}
|D^{\alpha}f(x)|+\|f\|_{\CKH}.
$$
Throughout the paper we let $S$ denote an arbitrary closed
subset of $\RN$.
\par In this paper we study the following extension
problem.
\par {\bf Problem.} Given a positive integer $k$ and an
arbitrary function $f:S\to\R$, what is a necessary and
sufficient condition for $f$ to be the restriction to $S$
of a function $F\in\CKO$? \vspace*{2mm}
\par This is a variant of a classical problem which is
known in the literature as the Whitney Extension Problem
\cite {W1,W2}. It has attracted a lot of attention in
recent years. We refer the reader to \cite{BS1}-\cite{BS4},
\cite{F1}-\cite{F7}, \cite{BMP1,BMP2} and \cite{Z1,Z2} and
references therein for numerous results in this direction,
and for a variety of techniques for obtaining them.
\par This note is devoted to the {\it phenomenon of
``finiteness"} in the Whitney problem. It turns out that,
in many cases, Whitney-type problems for different spaces
of smooth functions can be reduced to the same kinds of
problems, but for finite sets with {\it prescribed numbers
of points.}
\par For the space $\COM$ (with
$\omega(t)=t^p,\,0<p\le 1,$) and for the Zygmund space,
this phenomenon has been studied in the author's papers
\cite{S0,S1}. The case of an arbitrary $\omega$ was treated
in joint papers with Yu.\ Brudnyi \cite{BS1,BS4}. It was
shown that {\it a function $f$ defined on $S$ can be
extended to a function $F\in \COM$ with $\|F\|_{\COM}\le
\gamma=\gamma(n)$ provided its restriction $f|_{S'}$ to
every subset $S'\subset S$ consisting of at most
$N(n)=3\cdot 2^{n-1}$ points can be extended to a function
$F_{S'}\in \COM$ with $\|F_{S'}\|_{\COM}\le 1$.} (Moreover,
the value $3\cdot 2^{n-1}$ is sharp \cite{S1,BS4}.)
\par This result is an example of {\it ``the
finiteness property"} of the space $\COM$. We call the
number $N$ appearing in formulations of finiteness
properties {\it ``the finiteness number".}
\par In his pioneering work \cite{W2}, H. Whitney
characterized the restriction of the space $C^k(\R), k\ge
1,$ to an arbitrary subset $S\subset\R$ in terms of divided
differences of functions. An application of Whitney's
method to the space $\COR$ implies the finiteness property
for this space with the finiteness number $N=k+2$.
\par An impressive breakthrough in the solution of the Whitney
problem for $C^{k,\omega}$-spaces has recently been made by
C. Fefferman \cite{F1}-\cite{F7}. In this paper we will
consider two of his remarkable results related to the
finiteness property and its generalizations for the space
$\CKO$. Here is the first of them:
\begin{theorem}\label{FNT}(C. Fefferman \cite{F1,F3}).
There is a positive integer $N=N(k,n)$ such that the
following is true: Suppose we are given a function
$\omega$, a set $S\subset \RN$, and  functions
$f:S\rightarrow \R$ and $\xi:S\rightarrow \R_{+}$. Assume
that, for any $S'\subset S$ with at most $N$ points, there
exists a function $F_{S'}\in \CKO$ with
$\|F_{S'}\|_{\CKO}\le 1$, and
$$
|F_{S'}(x)-f(x)| \leq \xi(x)~~~for~all~~~x\in S'.
$$
Then there exists $F\in C^{k,\omega}(\RN)$, with
$\|F\|_{C^{k,\omega}(\RN)}\leq \gamma$ and
$$
|F(x)-f(x)| \leq \gamma\cdot\xi(x),~~~~x\in S.
$$
Here $\gamma=\gamma(k,n)$ is a constant depending only on
$k$ and $n$.
\end{theorem}
\par In particular, if the function $\xi $ is chosen to be
identically zero, Theorem \ref{FNT} shows that the space
$\CKO$ possesses the finiteness property {\it for all}
$k,n\ge 1$.
\par An upper bound for the finiteness number $N(k,n)$
given in \cite{F1,F3} is
\bel{Ne1}
N(k,n)\le(\dim\PK + 1)^ {3\cdot2^{\dim\PK}}.
\ee
Here $\PK$ stands for the space of polynomials of degree at
most $k$ defined on $\RN$. (Recall that
$\dim\PK={n+k\choose k}$.)
\par Our first result, Theorem \ref{TF3}, states
that the expression bounding $N(k,n)$ in \rf{Ne1} can be
replaced by a considerably smaller expression which depends
on $\dim\PK$ exponentially.
\begin{theorem}\label{TF3}
Theorem \ref{FNT} holds with the finiteness number
$N(k,n)=\,2^{\dim\PK}$.
\end{theorem}
\begin{remark} {\em  In the spring of 2005 I learned that
E. Bierstone and P. Milman  obtained an improvement of
estimate \rf{Ne1} to a bound that is exponential. Recently
P. Milman kindly drew my attention to the fact that their
result gives precisely the estimate $2^{\dim\PK}$ for the
spaces $C^k(\RN)$ and $\CKO$.}
\end{remark}
\par In fact there are many different versions of the
Whitney extension problem. These versions arise when one
considers a possibly different space of smooth functions on
$\RN$ and a possibly different collection of given
information about the function on the set $S$. In his
classical paper \cite{W1}, Whitney solved a version for the
space $C^k(\RN)$ in the case where the given information
about the function includes its values and the values of
{\it all} of its partial derivatives of all orders up to
$k$ on the set $S$. Using Whitney's extension method G.
Glaeser \cite{G} proved a similar result for the space
$\CKO$. Let us recall its formulation.
\par Given a $k$-times differentiable
function $f$ and $x\in\RN,$ we let $T_{x}^{k}(f)$ denote
the Taylor polynomial of $f$ at $x$ of degree at most $k$:
$$
T_{x}^{k}(f)(y):=\sum_{|\alpha|\leq
k}\frac{1}{\alpha!}(D^{\alpha}f)(x)(y-x)^{\alpha}~,~~y\in
\RN.
$$
\begin{theorem}\label{WG}(Whitney-Glaeser).
Given a family of polynomials
$\{P_x\in\PK:x\in S\}$
there is a function $F\in\CKO$ such that
$T_{x}^{k}(F)=P_{x}$ for every $x\in S$ if and only if
there is a constant $\lambda>0$ such that for every
$\alpha, |\alpha|\le k$ we have
\bel{nd} |D^\alpha P_x(x)|\le \lambda ~~~~for~all~~~x\in S,
\ee
and
\bel{dp}
\max\{|D^\alpha(P_x-P_y)(x)|,|D^\alpha(P_x-P_y)(y)|\} \le
\lambda\, \|x-y\|^{k-|\alpha|}\,\omega(\|x-y\|), \ee
for all $x,y\in S$. Moreover,
$$
\inf\{\|F\|_{\CKO}:~T_{x}^{k}(F)=P_{x} ,x\in S\}\approx
\inf\lambda
$$
with constants of equivalence depending only on $k$ and
$n$.
\end{theorem}
\par Observe that Theorem \ref{WG} can be interpreted as a
finiteness theorem with the finiteness number $N=2$. In
fact, the inequalities \rf{nd} and \rf{dp} depend on at
most $2$ (arbitrary) points of $S$ so that the sufficiency
part of this result can be reformulated as follows: There
is a function $F\in\CKO$ with $\|F\|_{\CKO}\le\gamma(k,n)$
satisfying $T_{x}^{k}(F)=P_{x},~x\in S,$ provided for every
two-point set $S'\subset S$ there exists a function
$F_{S'}\in\CKO$ with $\|F_{S'}\|_{\CKO}\le 1$ such that
$T_{x}^{k}(F_{S'})=P_{x} ,~x\in S'.$
\par In \cite{F4} C. Fefferman considered a version of the
Whitney problem in which the family of polynomials
$\{P_x\in\PK:x\in S\}$ is replaced by a family $\{G(x):x\in
S\}$ of convex centrally-symmetric subsets of $\PK$. He
raised the following question: How can we decide whether
there exist $F\in\CKO$ and a constant $A>0$ such that
$$
T_{x}^{k}(F)\in A\circledcirc G(x)
 ~~~~{\rm for~all}~~x\in S~?
$$
Here $A\circledcirc G(x)$ denotes the dilation of $G(x)$
with respect to its center by a factor of $A$.
\par Let $P_x\in\PK$ be the center of the set $G(x)$.
This means that $G(x)$ can be represented in the form
$G(x)=P_x+\sigma(x)$ where $\sigma(x)\subset\PK$ is a
convex family of polynomials which is centrally symmetric
with respect to $0$. It is shown in \cite{F4} that, under
certain conditions on the sets $\sigma(x)$, the finiteness
property holds. We say that a set $\sigma(x)\subset\PK$ is
``Whitney $\omega$-convex" (with Whitney constant $A$) at
$x\in\RN$ if the following two conditions are satisfied:
\par (i). $\sigma(x)$ is closed, convex and symmetric with
respect to $0$;
\par (ii). Suppose $P\in \sigma(x)$, $Q\in
\PK$ and $\delta \in (0,1]$. Assume that $P$ and $Q$
satisfy the estimates
$$ |\partial^{\beta}P(x)|\leq
\omega(\delta)\delta^{k-|\beta|} ~~~{\rm and} ~~~~
|\partial^{\beta}Q(x)|\leq \delta^{-|\beta|}
$$
for all $|\beta|\leq k.$ Then $T^k_x(P\cdot Q)\in
A\sigma(x)$. (See \cite{F4}, p. 579.)
\begin{theorem}\label{F4}(\cite{F4}) Given integers
$k,n\ge 1$ there is a constant $N=N(k,n)$ for which the
following holds: For each $x\in S$, suppose we are given a
polynomial $P_x\in\PK$, and a Whitney $\omega$-convex set
$\sigma(x)$ with Whitney constant $A$. Suppose that for
every subset $S'$ of $S$ with cardinality at most $N$ there
exists a function $F_{S'}\in \CKO$ such that
$\|F_{S'}\|_{\CKO}\leq 1$ and
$$
T_{x}^{k}(F_{S'})\in P_{x}+\sigma(x)~~~~for~all~~~x\in S'.
$$
Then there exists a function $F\in \CKO$, satisfying
$\|F\|_{\CKO}\leq \gamma$ and
\bel{dil}
 T_{x}^{k}(F)\in P_{x}+\gamma\cdot\sigma(x)~,~~~~x\in S.
\ee
Here $\gamma$ depends only on $k,n$ and the Whitney
constant $A$.
\end{theorem}
\par A particular case of this result for
$\sigma(x)=\{P\in\PK:D^{\alpha}P(x)=0,\,|\alpha|\le k-1\}$
and $\omega(t)=t^p,\,0<p\le 1,$ with the finiteness number
$N=3\cdot 2^{{n+k-2\choose k}}$ has been proved in
\cite{BS3}.
\par Analogously to Theorem \ref{TF3}, our second result
in this paper gives an explicit upper bound for a
finiteness number.
\begin{theorem}\label{TF5}
Theorem \ref{F4} holds with the finiteness number
$$N(k,n)=\,2^{\min\{\ell+1,\,\dim\PK\}},$$
where $\ell=\max_{x\in S}{\dim \sigma(x)}$.
\end{theorem}
\begin{remark}
{\em Of course $\ell$ necessarily satisfies $\ell\le
\dim\PK={n+k\choose k}$.}
\end{remark}
\par In fact both of our new estimates for finiteness
numbers are corollaries of the following theorem which is
the main result of this paper.
\begin{theorem}\label{TF1} Let $G$ be a mapping
defined on a finite set $S\subset\RN$ which assigns a
convex set of polynomials $G(x)\subset\PK$ of dimension at
most $\ell$ to every point $x$ of $S$. Suppose that, for
every subset $S'$ of $S$ consisting of at most
$2^{\,\min\{\ell+1,\,\dim\PK\}}$ points, there exists a
function $F_{S'}\in \CKO$ such that $\|F_{S'}\|_{\CKO}\leq
1$ and $T_{x}^{k}(F_{S'})\in G(x)$ for all $x\in S'$. Then
there is a function $F\in\CKO$, satisfying
$\|F\|_{\CKO}\le\gamma$ and
$$
T_{x}^{k}(F)\in G(x) ~~~for~all~~~x\in S.
$$
Here $\gamma$ depends only on $k,n$ and $\card S$.
\end{theorem}
\par Comparing this result with Theorem
\ref{F4}, let us note that here there are no restrictions
on $G$. Moreover, here $T_{x}^{k}(F)$ belongs to $G(x)$
itself and not merely to its dilation as in \rf{dil}.
However the price of that we have to pay to obtain such a
general result is that we have to permit the constant
$\gamma$ (controlling the $C^{k,\omega}$-norm of the
function $F$) to {\it depend on the number of points of
$S$.}
\par We can use the rather informal and imprecise terminology
``$\CKO$ has the weak finiteness property" to express the
kind of result obtained in Theorem \ref{TF1} where $\gamma$
depends on the number of points of $S$. The fact that such
a weak finiteness property holds, strongly suggests that we
can reasonably hope to establish an analogous ``strong
finiteness property", by which we mean a result with
$\gamma$ depending only on $k$ and $n$. Such a result may
possibly require some additional very mild conditions to be
imposed on the mapping $G$.
\par The weak finiteness property also provides an upper
bound for the finiteness constant whenever the strong
finiteness property holds. For instance, Fefferman's
Theorems \ref{FNT} and \ref{F4} reduce the problem to a set
of cardinality at most $N(k,n)$ while the weak finiteness
property decreases this number to $2^{\dim\PK}$ (as in
Theorem \ref{TF3}) or to $2^{\min\{l+1,\dim\PK\}}$ (Theorem
\ref{TF5}).
\par We prove Theorem \ref{TF1} in Section 4. The proof is
based on an approach presented in Sections 2 and 3.
\par The crucial ingredient in this approach
is an isomorphism between the space $\CKO|_S$ and a certain
space of Lipschitz mappings from $S$ into the product
$\PRN$ equipped with a certain metric $\dm$. We define
$\dm$ and study its properties in Section 2. One of these
properties, which is obtained in Proposition \ref{dmcalc},
is a useful formula for calculating $\dm$, namely
$$
\dm(T,T')\approx \max_{|\alpha|\le
k}\left\{\omega(\|x-x'\|),\vfa(|D^\alpha
(P-P')(x)|),\vfa(|D^\alpha (P-P')(x')|)\right\}
$$
where $T=(P,x)$ and $T'=(P',x')$ are any two elements of
$\PRN$, and $\vfa:=\omega((s^{k-|\alpha|}\omega(s))^{-1})$.
We refer to the set
$$
\PRN:=\{(P,x):~P\in\PK,x\in\RN\}
$$
as {\it the space of (potential) $k$-jets.} This name and
also the definition of $\dm$ are motivated by the
Whitney-Glaeser extension theorem \ref{WG}.
\par  Given $T=(P,x)\in\PRN$ and $\lambda\in\R$ we define
$\lambda\circ T:=(\lambda P,x)$. Then inequality \rf{dp} of
the Whitney-Glaeser extension theorem can be reformulated
as follows:
\bel{M1}
\dm(\lambda^{-1}\circ T_x,\lambda^{-1}\circ
T_y)\le \omega(\|x-y\|),
~~~T_x:=(P_x,x),T_y:=(P_y,y)\in\PK\times S.
\ee
\par We define a metric on $S$ by setting
$\ro(x,y):=\omega(\|x-y\|)$ for all $x,y\in S$ and we let
$\SO$ be the metric space $\SO:=(S,\ro)$. We also consider
$\PRN$ as a metric space with respect to $\dm$, i.e., we
set $\TK:=(\PRN,\dm)$. Let $\Lip(\SO,\TK)$ denote the space
of Lipschitz mappings from $S$ (equipped with the metric
$\ro$) into $\PRN$ (with the metric $\dm$). Inequality
\rf{M1} motivates us to equip this space with a ``norm" by
setting
\bel{DO}
\|T\|_{LO(S)}:=\inf\{\lambda:~\|\lambda^{-1}\circ
T\|_{\Lip(\SO,\TK)}\le 1\}.
\ee
\par We call $\|\cdot\|_{LO(S)}$ the Lipschitz-Orlicz norm.
We use it to define a second ``norm" by setting
\bel{TS} \|T\|^*_{LO(S)}:=\max_{|\alpha|\le k}\,\sup_{x\in
S} |D^\alpha P_x(x)|+\|T\|_{LO(S)} \ee
and we introduce the subspace $\BLip(\SO,\TK)$ of
$\Lip(\SO,\TK)$ of ``bounded" Lipschitz mappings
$T(x)=(P_x,z_x),~x\in S,$ defined by the finiteness of the
``norm" \rf{TS}.
\par Now the Whitney-Glaeser extension theorem implies
the following
\begin{proposition}\label{WLip}
Given a family of polynomials $\{P_x\in\PK:~x\in S\}$,
there is a function $F\in\CKO$ such that $T^k_x(F)=P_x$ for
every $x\in S$ if and only if the mapping
$T(x):=(P_x,x),x\in S,$ belongs to $\BLip(\SO,\TK)$.
Moreover,
$$
\inf\{\|F\|_{\CKO}:~T^k_x(F)=P_x,x\in
S\}\approx\|T\|^*_{LO(S)}
$$
with constants of equivalence depending only on $k$ and
$n$.
\end{proposition}
\par Applying this proposition to $S=\RN$ we obtain an
interesting isomorphism between $\CKO$ and a certain
subfamily of $\Lip(\RN_\omega,\TK)$. Namely, every function
$F\in\CKO$ gives rise to a Lipschitz mapping from
$\RN_\omega:=(\RN,\omega(\|\cdot\|))$ into $\PRN$ defined
by the formula $ T(x):=(T^k_x(F),x), x\in\RN.$ On the other
hand, every Lipschitz mapping from $\RN_\omega$ to $\PRN$
of the form $T(x):=(P_x,x), x\in\RN,$ generates a function
$F(x):=P_x(x), x\in\RN,$ such that $F\in\CKO$ and
$T^k_x(F)=P_x, x\in\RN$.
\par Let us restate this more concisely: The mapping
$$
\CKO\ni F \mapsto T(x):=(T^k_x(F),x)\in
\Lip(\RN_\omega,\TK)
$$
and its inverse mapping
$$
\Lip(\RN_\omega,\TK)\ni T(x):=(P_x,x)\mapsto
F(x):=P_x(x)\in\CKO
$$
provide an isomorphism between $\CKO$ and the subfamily of
$\Lip(\RN_\omega,\TK)$ consisting of all elements of the
form $T(x):=(P_x,x),~x\in\RN$. Moreover, Proposition
\ref{WLip} states that this isomorphism in some sense
``preserves restrictions".
\par The above ideas and results are presented in Section 2.
They show that even though Whitney's problem deals with
restrictions of $k$-times differentiable functions, it is
also a problem about Lipschitz mappings defined on subsets
of $\RN$ and taking values in a very non-linear metric
space $\TK=(\PRN,\dm)$. More specifically, the Whitney
problem can be reformulated as a problem about {\it
Lipschitz selections of set-valued mappings} from $S$ into
$2^{\TK}$. We study this problem in Section 3. We remark
that the Lipschitz selection method has already been used
to obtain a solution  to the Whitney problem for the space
$\COM$, see \cite{S1,S3, BS4}.
\par We recall some relevant definitions: Let
$X=(\Mc,\rho)$ and $Y=(\Tc,d)$ be metric spaces and let
$\Gc:\Mc \to 2^{\Tc}$ be a set-valued mapping, i.e., a
mapping which assigns a {\it subset} $\Gc(x)\subset\Tc$ to
each $x\in\Mc$. A function $g:\Mc\to\Tc$ is said to be a
{\it selection} of $\Gc$ if $g(x)\in \Gc(x)$ for all
$x\in\Mc$. If a selection $g$ is an element of $\Lip(X,Y)$
then it is said to be a {\it Lipschitz selection} of the
mapping $\Gc$.  (For various results and techniques related
to the problem of the existence of Lipschitz selections in
the case where $Y=(\Tc,d)$ is a Banach space, we refer the
reader to \cite{S2,S3,S4} and references therein.)
\par It turns out that Theorem \ref{TF1},
the ``weak finiteness" theorem, is equivalent to the
following Helly-type criterion for the existence of a
Lipschitz selection.
\begin{theorem}\label{FW} Let $S\subset\RN$ be a finite set
and let $\Gc(x)=(G(x),x),x\in S,$ be a set-valued mapping
such that for each $x\in S$ the set $G(x)\subset\PK$ is a
convex set of polynomials of dimension at most $\ell$.
Suppose that there exists a constant $K>0$ such that, for
every subset $S'\subset S$ consisting of at most
$2^{\min\{\ell+1,\dim\PK\}}$ points, the restriction
$\Gc|_{S'}$ has a Lipschitz selection
$g_{S'}\in\Lip(S',\TK)$ with $\|g_{S'}\|_{LO(S')}\le K$.
Then $\Gc$, considered as a map on all of $S$, has a
Lipschitz selection $g\in\Lip(S,\TK)$ with
$\|g\|_{LO(S)}\le \gamma K$, where the constant $\gamma$
depends only on $k,n$ and $\card S$.
\end{theorem}
\par The proof of this result relies on some methods and
ideas developed for the case of set-valued mappings which
take their values in Banach spaces, see, e.g. Shvartsman
\cite{S2,S3,S4}. In particular, an analog of Theorem
\ref{FW} for Banach spaces has been proved in \cite{S3}.
Our strategy will be to adapt that proof to the case of the
metric space $\TK=(\PRN,\dm)$. As in the case of Banach
spaces  our adapted proof will be based on Helly's
intersection theorem \cite{DGK} and a combinatorial result
about a structure of finite metric graphs (Proposition
\ref{MT}).
\par {\bf Acknowledgment.} I am greatly indebted to Michael
Cwikel, Charles Fefferman and Naum Zobin  for interesting
discussions and helpful suggestions and remarks.
\SECT{2. $\CKO$ as a space of Lipschitz mappings}{2}
\indent
\par The point of departure for our approach is inequality
\rf{dp} of the Whitney extension theorem. This inequality
motivates the definition of a certain special metric on the
set $\PRN$ which allows us to identify the restriction
$\CKO|_S$ with a space of {\it Lipschitz mappings} from $S$
into $\PRN$.
\par Observe that without loss of generality we may assume
that $\omega$ is a {\it strictly increasing} concave
function on $\R_+$. (In fact, for every positive concave
$\omega:\R_+\to\R_+$ there is a concave strictly increasing
function $\omega^*$ such that $\omega^*\le\omega$.
Therefore $\tilde{\omega}:=\omega+\omega^*$ is a concave
strictly increasing function satisfying
$\omega\le\tilde{\omega}\le 2\omega$.)
\par Now let us define a metric on $\PRN$. To this end given
multiindex $\alpha,|\alpha|\le k,$ we define a function
$\vfa:\R_+\to\R_+$ by letting
\bel{vfa}
\vfa:=\omega((s^{k-|\alpha|}\omega(s))^{-1})
\ee
for $|\alpha|<k$ and $\vfa(t):=t$ for $|\alpha|=k$. Since
for $|\alpha|<k$ the function $s^{k-|\alpha|}\omega(s)$ is
strictly increasing, the inverse function
$\psi_{\alpha}:=(s^{k-|\alpha|}\omega(s))^{-1}$ is
well-defined so that the function
$\vfa=\omega(\psi_{\alpha})$ is well-defined as well. It
can be also readily seen that
$$
\vfa'(t)=\frac{1}{\psi^{k-|\alpha|}_{\alpha}(t)
+(k-|\alpha|)\psi_{\alpha}^{k-|\alpha|-1}(t)
\omega(\psi_{\alpha}(t))/\omega'(\psi_{\alpha}(t))}.
$$
Since $\psi_{\alpha}$ and $\omega$ are non-decreasing and
$\omega'$ is non-increasing, $\vfa'$ is non-increasing, so
that $\vfa$ is a {\it concave} function.
\par Fix two $k$-jets $T_0=(P_0,x_0),
T_1=(P_1, x_1)\in\PRN$ and put
\bel{deffi}
\dom(T_0,T_1):=\max_{|\alpha|\le k}
\{\omega(\|x_0-x_1\|), \vfa(|D^\alpha(P_0-P_1)(x_0)|),
\vfa(|D^\alpha(P_0-P_1)(x_1)|\}.
\ee
Clearly,
\bel{pom}
\dom((P,x_0),(P,x_1))=\omega(\|x_0-x_1\|) ~~~{\rm
for~every}~~~P\in\PK,~x,y\in\RN.
\ee
\par Recall that $\lambda\circ T:=(\lambda P,x)$ where
$T=(P,x)\in\PRN$ and $\lambda\in\R$. In these settings
inequality \rf{dp} of the Whitney-Glaeser extension theorem
means the following
$$
\dom(\lambda^{-1}\circ T_x,\lambda^{-1}\circ T_y)\le
\omega(\|x-y\|), ~~~~~~T_x:=(P_x,x),T_y:=(P_y,y).
$$
\par Now we define a metric $\dm$ on $\PRN$ by letting
\bel{dm} \dm(T,T'):=\inf\sum_{i=0}^{m-1}\dom(T_i,T_{i+1})
\ee
where the infimum is taken over all finite families
$\{T_0,T_1,...,T_m\}\subset\PRN$ such that $T_0=T$ and
$T_m=T'$.
\par In particular, since $\omega$ is subadditive, by this
definition for every $x,y\in\RN$
$$
\dm((P,x),(Q,y))\ge \omega(\|x-y\|), ~~~~P,Q\in\PK,
$$
and  by \rf{pom}
\bel{pdm}
\dm((P,x),(P,y))=\omega(\|x-y\|) ~~~{\rm
for~every}~~~P\in\PK.
\ee
\par The main result of this section is the following
\begin{theorem}\label{DOC}
For every $T,T'\in\PRN$ we have
$$
\dm(T,T')\le \dom(T,T')\le \dm(e^n\circ T,e^n\circ T').
$$
\end{theorem}
\par A proof of the theorem is based on a series of
auxiliary lemmas.
\begin{lemma}\label{fiab}
For every $t_1,t_2>0$ and every multiindexes $\alpha,
\beta$ such that $|\alpha|+|\beta|\le k$ we have
$$
\vfa(t_1^{|\beta|}t_2)\le \max\{\omega(t_1),
\vf_{\alpha+\beta}(t_2)\}.
$$
\end{lemma}
\par{\it Proof.} If $|\alpha|=k$, then $\beta=0$ so that
nothing to prove. Therefore we will assume that
$|\alpha|<k$. In this case by \rf{vfa}
 $\vfa(t_1^{|\beta|}t_2) =\omega(u)$ where
$u:=(s^{k-|\alpha|}\omega(s))^{-1}(t_1^{|\beta|}t_2).$
Hence $t_1^{|\beta|}t_2=u^{k-|\alpha|}\omega(u).$
\par Suppose that $|\alpha|+|\beta|< k$. Then
$\vf_{\alpha+\beta}(t_2)=\omega(v)$ where
$v:=(s^{k-|\alpha|-|\beta|}\omega(s))^{-1}(t_2)$ so that
$t_2=v^{k-|\alpha|-|\beta|}\omega(v).$
\par Put $w:=\max\{t_1,v\}.$ Then
\bel{max} \max\{\omega(t_1),\vf_{\alpha+\beta}(t_2)\}
=\max\{\omega(t_1),\omega(v)\}=
\omega(\max\{t_1,v\})=\omega(w). \ee
Hence
$$
t_1^{|\beta|}t_2=t_1^{|\beta|}
v^{k-|\alpha|-|\beta|}\omega(v) \le w^{|\beta|}
w^{k-|\alpha|-|\beta|}\omega(w) =w^{k-|\alpha|}\omega(w).
$$
But $t_1^{|\beta|}t_2=u^{k-|\alpha|}\omega(u)$ so that
$$
u^{k-|\alpha|}\omega(u)=t_1^{|\beta|}t_2 \le
w^{k-|\alpha|}\omega(w).
$$
Since $t^{k-|\alpha|}\omega(t)$ is an strictly increasing
function, this implies $u\le w$. Recall also that
$\vfa(t_1^{|\beta|}t_2)=\omega(u)$. Then by \rf{max}
$$
\vfa(t_1^{|\beta|}t_2)=\omega(u)\le \omega(w)
=\max\{\omega(t_1),\vf_{\alpha+\beta}(t_2)\}.
$$
\par It remains to consider the case $|\alpha|+|\beta|=k$.
In this case by the definition of $\vfa$ we have
$\vf_{\alpha+\beta}(t_2)=t_2.$ If $\sup_{t>0}\omega(t)\le
t_2$, then
$$
\vfa(t_1^{|\beta|}t_2)=\omega(u)\le
t_2=\vf_{\alpha+\beta}(t_2)
$$
and the lemma follows. If $t_2<\sup_{t>0}\omega(t)$ then
there is $v>0$ such that $t_2=\omega(v)$. This shows that
equality $t_2=v^{k-|\alpha|-|\beta|}\omega(v)$ holds for
the case $|\alpha|+|\beta|=k$ as well. The lemma is
proved.\bx
\begin{lemma}\label{p-norm}
Let $Q\in\PK$ and let $a,b\in\RN$. Then for every
multiindex $\alpha, |\alpha|\le k,$ we have
$$
|D^\alpha Q(b)|\le \sum_{|\beta|\le k-|\alpha|}
\frac{1}{\beta!}|D^{\alpha+\beta}Q(a)|\cdot
\|b-a\|^{|\beta|}.
$$
\end{lemma}
\par {\it Proof.} Since $D^\alpha Q$ is a polynomial of
degree at most $k-|\alpha|$, Taylor's formula for $D^\alpha
Q$ at $a$ gives
$$
D^\alpha Q(b)= \sum_{|\beta|\le k-|\alpha|}
\frac{1}{\beta!}D^{\alpha+\beta} Q(a)(b-a)^\beta
$$
which immediately implies the required inequality of the
lemma.\bx
\begin{lemma}\label{mainest} Let
$T=(P,x),T'=(P',x')\in\PRN$ and let
$$\{T_i=(P_i,x_i):~i=0,1,...,m\}$$
be a finite family of elements of $\PRN$ such that
$T_0=T,\,T_m=T'$. Then for every $\alpha, |\alpha|\le k,$
$$
|D^\alpha(P-P')(x)|\le e^n\max_{|\beta|\le k-|\alpha|}
\sum^{m-1}_{i=0}|D^{\alpha+\beta}
(P_i-P_{i+1})(x_i)|\cdot\|x-x_i\|^{|\beta|}.
$$
\end{lemma}
\par {\it Proof.}
By Lemma \ref{p-norm}
$$
|D^\alpha (P_i-P_{i+1})(x)|\le \sum_{|\beta|\le k-|\alpha|}
\frac{1}{\beta!}|D^{\alpha+\beta}(P_i-P_{i+1})(x_i)|\cdot
\|x-x_i\|^{|\beta|}
$$
so that
\be |D^\alpha(P-P')(x)|&\le& \sum^{m-1}_{i=0}|D^\alpha
(P_i-P_{i+1})(x)|
\nn\\
&\le& \sum^{m-1}_{i=0}\sum_{|\beta|\le
k-|\alpha|}\frac{1}{\beta!}|D^{\alpha+\beta}
(P_i-P_{i+1})(x_i)|\cdot\|x-x_i\|^{|\beta|} \nn\\
&=& \sum_{|\beta|\le k-|\alpha|}\frac{1}{\beta!}
\sum^{m-1}_{i=0}|D^{\alpha+\beta}
(P_i-P_{i+1})(x_i)|\cdot\|x-x_i\|^{|\beta|}. \nn\ee
Hence
\be |D^\alpha(P-P')(x)| &\le& \left(\sum_{|\beta|\le
k-|\alpha|}\frac{1}{\beta!}\right) \max_{|\beta|\le
k-|\alpha|} \sum^{m-1}_{i=0}|D^{\alpha+\beta}
(P_i-P_{i+1})(x_i)|\cdot\|x-x_i\|^{|\beta|}\nn\\
&\le&
e^n\max_{|\beta|\le k-|\alpha|}
\sum^{m-1}_{i=0}|D^{\alpha+\beta}
(P_i-P_{i+1})(x_i)|\cdot\|x-x_i\|^{|\beta|}
 \nn \ee
proving the lemma.\bx
\par We are in a position to prove Theorem
\ref{DOC}.\medskip
\par {\it Proof of Theorem \ref{DOC}.} The first inequality
follows from definition \rf{dm}. Let us prove the second
inequality. Consider a family
$\{T_i=(P_i,x_i):~i=0,1,...,m\}\subset\PRN$ such that
$T_0=T:=(P,x),T_m=T':=(P',x')$. Thus $P_i\in\PK,~x_i\in\RN$
for every $i=0,...,m$ and $P_0=P, P_m=P', x_0=x, x_m=x'$.
Let us prove that
\bel{dt}
\dom(T,T')\le \sum^{m-1}_{i=0}\dom(e^n\circ
T_i,e^n\circ T_{i+1}).
\ee
\par Let us fix a multiindex $\alpha,~|\alpha|\le k$ and
estimate $\vfa(|D^\alpha(P-P')(x)|)$.
By Lemma \ref{mainest}
$$
|D^\alpha(P-P')(x)|\le e^n\max_{|\beta|\le k-|\alpha|}
\sum^{m-1}_{i=0}|D^{\alpha+\beta}
(P_i-P_{i+1})(x_i)|\cdot\|x-x_i\|^{|\beta|}.
$$
Put $\tP_i:=e^n P_i$. Then the latter inequality implies
$$
|D^\alpha(P-P')(x)|\le \max_{|\beta|\le k-|\alpha|}
\sum^{m-1}_{i=0}|D^{\alpha+\beta}
(\tP_i-\tP_{i+1})(x_i)|\cdot\|x-x_i\|^{|\beta|}.
$$
Since
$\|x-x_i\|\le\sum^{m-1}_{i=0} \|x_i-x_{i+1}\| $
we obtain
$$
|D^\alpha(P-P')(x)|\le \max_{|\beta|\le k-|\alpha|}
\left(\sum^{m-1}_{i=0}\|x_i-x_{i+1}\|\right)^{|\beta|}\cdot
\left(\sum^{m-1}_{i=0}|D^{\alpha+\beta}(\tP_i-
\tP_{i+1})(x_i)|\right).
$$
Recall that the function $\vfa$ defined by \rf{vfa} is
non-decreasing. Hence
\be
&&\vfa(|D^\alpha(P-P')(x)|)\nn\\
&&\le\max_{|\beta|\le k-|\alpha|}
\vfa\left(\left(\sum^{m-1}_{i=0}\|x_i-x_{i+1}\|
\right)^{|\beta|}\cdot
\left(\sum^{m-1}_{i=0}
|D^{\alpha+\beta}(\tP_i-\tP_{i+1})(x_i)|\right)\right).
\nn \ee
By Lemma \ref{fiab}
\be
&&\vfa\left(\left(\sum^{m-1}_{i=0}
\|x_i-x_{i+1}\|\right)^{|\beta|}\cdot
\left(\sum^{m-1}_{i=0}|D^{\alpha+\beta}
(\tP_i-\tP_{i+1})(x_i)|\right)\right)
\nn\\
&&\le
\max\left\{\omega\left(\sum^{m-1}_{i=0}
\|x_i-x_{i+1}\|\right),
\vf_{\alpha+\beta} \left(\sum^{m-1}_{i=0}|D^{\alpha+\beta}
(\tP_i-\tP_{i+1})(x_i)|\right)\right\}. \nn \ee
Since $\omega$ and $\vf_{\alpha+\beta}$ are  concave
functions on $\R_+$, they are subadditive so that
\bel{som}
\omega\left(\sum^{m-1}_{i=0}\|x_i-x_{i+1}\|\right) \le
\sum^{m-1}_{i=0}\omega(\|x_i-x_{i+1}\|) \ee
and
$$
\vf_{\alpha+\beta}
\left(\sum^{m-1}_{i=0}|D^{\alpha+\beta}
(\tP_i-\tP_{i+1})(x_i)|\right)
\le \sum^{m-1}_{i=0}\vf_{\alpha+\beta}
(|D^{\alpha+\beta}(\tP_i-\tP_{i+1})(x_i)|).
$$
Hence by definition \rf{deffi} of $\dom$ we have
\be
&&\vfa(|D^\alpha(P-P')(x)|)\nn\\
&&\le\max_{|\beta|\le k-|\alpha|}
 \left\{\sum^{m-1}_{i=0}\omega(\|x_i-x_{i+1}\|),
\sum^{m-1}_{i=0}\vf_{\alpha+\beta}
(|D^{\alpha+\beta}(\tP_i-\tP_{i+1})(x_i)|)\right\}\nn\\
&&\le\max_{|\beta|\le k-|\alpha|}
 \left\{\sum^{m-1}_{i=0}\max\{\omega(\|x_i-x_{i+1}\|),
\vf_{\alpha+\beta}
(|D^{\alpha+\beta}(\tP_i-\tP_{i+1})(x_i)|)\}\right\}\nn\\
&&\le
\sum^{m-1}_{i=0}\max_{|\beta|\le k-|\alpha|}
 \{\omega(\|x_i-x_{i+1}\|),
\vf_{\alpha+\beta}
(|D^{\alpha+\beta}(\tP_i-\tP_{i+1})(x_i)|)\} \nn\\
&&\le
\sum^{m-1}_{i=0}\dom((\tP_i,x_i),(\tP_{i+1},x_{i+1}))=
\sum^{m-1}_{i=0}\dom(e^n\circ T_i,e^n\circ T_{i+1}). \nn
\ee
(Recall that $\tP_i:=e^n P_i$).
In a similar way we prove that
$$
\vfa(|D^\alpha(P-P')(x')|) \le
\sum^{m-1}_{i=0}\dom(e^n\circ T_i,e^n\circ T_{i+1}).
$$
Combining this with \rf{som} and definition \rf{deffi} we
obtain the required inequality \rf{dt}. By this inequality
$$
\dom(T,T')\le \inf \sum_{i=0}^{m-1} \dom(e^n\circ T_i,
e^n\circ T_{i+1})
$$
where the infimum is taken over all families
$\{T_0,T_1,...,T_m\}\subset\PRN$ such that $T_0=T$ and
$T_m=T'$. By \rf{dm} this infimum is equal to
 $\dm(e^n\circ T,e^n\circ T')$ and the proof
is finished.\bx
\par The following proposition presents two formulae for
calculation of the metric $\dm$.
\begin{proposition}\label{dmcalc}
For every $T=(P,x),T'=(P',x')\in\PRN$ we have
\par (i).
$$
e^{-n}\dom(T,T')\le\dm(T,T')\le\dom(T,T');
$$
\par (ii).
$$
e^{-n}\dm(T,T')\le\max_{|\alpha|\le k} \{\omega(\|x-x'\|),
\vfa(|D^\alpha (P-P')(x)|)\}\le e^n\dm(T,T').
$$
\end{proposition}
\par {\it Proof.} (i). By Theorem \ref{DOC} for every
$T,T'\in\PRN$ we have
$$
\dm(T,T')\le \dom(T,T')\le \dm(e^n\circ T,e^n\circ T')
$$
so that $\dom(e^{-n}\circ T,e^{-n}\circ T')\le \dm(T,T').$
On the other hand, since $\vfa$ is a concave function, for
every $\lambda\ge 1$ we obtain
\be &&\dom(\lambda\circ T,\lambda\circ T')\nn\\&&:=
\max_{|\alpha|\le k} \{\omega(\|x-x'\|),
\vfa(|D^\alpha(\lambda P-\lambda P')(x)|),
\vfa(|D^\alpha(\lambda P-\lambda P')(x')|\}\nn\\
&&\le \lambda\max_{|\alpha|\le k} \{\omega(\|x-x'\|),
\vfa(|D^\alpha(P-P')(x)|),
\vfa(|D^\alpha(P-P')(x')|\}\nn\\
&&= \lambda\dom(T,T')\nn \ee
so that
$$
\dom(T,T')=\dom(e^n\circ(e^{-n}\circ
T),e^n\circ(e^{-n}\circ T'))\le e^n\dom(e^{-n}\circ
T,e^{-n}\circ T')\le e^n\dm(T,T')
$$
proving (i).
\par (ii). By (i) and \rf{deffi} we have to prove that
$\dom(T,T')\le e^n I$ where
$$
I:=\max_{|\alpha|\le k} \{\omega(\|x-x'\|), \vfa(|D^\alpha
(P-P')(x)|)\}.
$$
This is equivalent to the inequality
\bel{n2}
\vfa(|D^\alpha (P-P')(x')|)\le e^nI,~~~
|\alpha|\le k.
\ee
By Lemma \ref{p-norm}
\be |D^\alpha (P-P')(x')|&\le& \left(\sum_{|\beta|\le
k-|\alpha|}\frac{1}{\beta!}\right) \max_{|\beta|\le
k-|\alpha|}
|D^{\alpha+\beta}(P-P')(x)|\cdot\|x-x'\|^{|\beta|}\nn\\
&\le& e^n\max_{|\beta|\le k-|\alpha|}
|D^{\alpha+\beta}(P-P')(x)|\cdot\|x-x'\|^{|\beta|}.\nn
\ee
Since $\vfa(\lambda t)\le \lambda \vfa(t),\,\lambda\ge 1,$
this implies
$$
\vfa(|D^\alpha (P-P')(x')|)\le e^n \max_{|\beta|\le
k-|\alpha|} \vfa(|D^{\alpha+\beta}(P-P')(x)
|\cdot\|x-x'\|^{|\beta|}).
$$
But by Lemma \ref{fiab} for every $\beta,\, |\beta|\le
k-|\alpha|,$ we have
$$
\vfa(|D^{\alpha+\beta}(P-P')(x)|\,\|x-x'\|^{|\beta|}) \le
\max\{\omega(\|x-x'\|),
\vf_{\alpha+\beta}(|D^{\alpha+\beta}(P-P')(x)|)\} \le I
$$
proving \rf{n2} and the lemma.\bx
\par In the next section we will need the following variant
of the triangle inequality for $\dm$.
\begin{lemma}\label{seqest} Let
$\{T_i=(P_i,x_i):~i=0,1,...,m\}$ be a family of elements of
$\PRN$ such that
\bel{IE}
\dm(T_i,T_{i+1})\le\omega(\|x_i-x_{i+1}\|),~~i=0,...,m-1.
\ee
Suppose that for some $\lambda\ge 1$ we have
$$
\sum_{i=0}^{m-1}\|x_i-x_{i+1}\| \le\lambda\|x_0-x_{m}\|,
~~~and~~~\sum_{i=0}^{m-1}\omega(\|x_i-x_{i+1}\|)
\le\lambda\omega(\|x_0-x_{m}\|).
$$
Then
\bel{cl} \dm(\tau^{-1}\circ T_0,\tau^{-1}\circ
T_{m})\le\omega(\|x_0-x_{m}\|) \ee
where $\tau:=e^{2n}\lambda^{k+1}$.
\end{lemma}
\par {\it Proof.} By Lemma \ref{mainest}
$$
|D^\alpha(P_0-P_m)(x_0)|\le e^n\max_{|\beta|\le k-|\alpha|}
\sum^{m-1}_{i=0}|D^{\alpha+\beta}
(P_i-P_{i+1})(x_i)|\cdot\|x_0-x_i\|^{|\beta|}.
$$
By Theorem \ref{DOC} and \rf{IE}
$$
\dom(e^{-n}\circ T_i,e^{-n}\circ T_{i+1})\le
\dm(T_i,T_{i+1})\le\omega(\|x_i-x_{i+1}\|),~~i=0,...,m-1,
$$
so that by definition \rf{deffi} of $\dom$
$$
\varphi_{\alpha+\beta}(e^{-n}|D^{\alpha+\beta}
(P_i-P_{i+1})(x_i)|)\le\omega(\|x_i-x_{i+1}\|).
$$
Since $\omega$ is strictly increasing, by definition
\rf{vfa} of $\varphi_\alpha$ we have
$$
|D^{\alpha+\beta} (P_i-P_{i+1})(x_i)|\le
e^n\|x_i-x_{i+1}\|^{k-|\alpha|-|\beta|}
\omega(\|x_i-x_{i+1}\|).
$$
Hence
$$
|D^\alpha(P_0-P_m)(x_0)|\le e^{2n}\max_{|\beta|\le
k-|\alpha|}
\sum_{i=0}^{m-1}\|x_i-x_{i+1}\|^{k-|\alpha|-|\beta|}
\omega(\|x_i-x_{i+1}\|)\|x_0-x_i\|^{|\beta|}.
$$
Clearly, $\|x_l-x_j\|\le\sum_{i=0}^{m-1}\|x_i-x_{i+1}\|$
for every non-negative integers $l$ and $j$ so that
\be
&&|D^\alpha(P_0-P_m)(x_0)|\nn\\
&\le& e^{2n}\max_{|\beta|\le k-|\alpha|} \sum_{i=0}^{m-1}
\left(\sum_{i=0}^{m-1}\|x_i-x_{i+1}\|\right)
^{k-|\alpha|-|\beta|}
\left(\sum_{i=0}^{m-1}\|x_i-x_{i+1}\|\right)^{|\beta|}
\omega(\|x_i-x_{i+1}\|)\nn\\
&=& e^{2n}\left(\sum_{i=0}^{m-1}
\|x_i-x_{i+1}\|\right)^{k-|\alpha|}
\left(\sum_{i=0}^{m-1}\omega(\|x_i-x_{i+1}\|)\right)
\nn\\
&\le&
e^{2n}\lambda^{k+1-|\alpha|}\|x_0-x_{m}\|^{k-|\alpha|}
\omega(\|x_0-x_{m}\|) \le
e^{2n}\lambda^{k+1}\|x_0-x_{m}\|^{k-|\alpha|}
\omega(\|x_0-x_{m}\|). \nn\ee
In a similar way  we obtain
$$
|D^\alpha(P_0-P_m)(x_m)| \le
e^{2n}\lambda^{k+1}\|x_0-x_{m}\|^{k-|\alpha|}
\omega(\|x_0-x_{m}\|),~~~|\alpha|\le k.
$$
In view of definition \rf{deffi} this implies
$$
\dom((e^{2n}\lambda^{k+1})^{-1}\circ T_0,
(e^{2n}\lambda^{k+1})^{-1}\circ T_{m})
\le\omega(\|x_0-x_{m}\|).
$$
It remains to note that by Theorem \ref{DOC} $\dm\le\dom$
and the lemma follows.\bx
\par We turn to the proof of Proposition \ref{WLip}. As
usual given metric spaces $X=\MR$ and $Y=(\Tc,\dd)$ we let
$\Lip(X,Y)$ denote the space of Lipschitz mappings from
$\Mc$ into $\Tc$. This (in general non-linear) space of
mappings $F:\Mc\to\Tc$ is equipped with the standard
``seminorm"
$$
\|F\|_{\Lip(X,Y)}:=\inf\{\lambda:~\dd(F(x),F(y))
\le\lambda\rho(x,y),~x,y\in\Mc\}.
$$
\par Recall that $\TK:=(\PRN,\dm)$ and $\SO:=(S,\ro)$ where
$\ro(x,y):=\omega(\|x-y\|)$, $x,y\in S$. Recall also  that
the space $\Lip(\SO,\TK)$ is normalized by the
Lipschitz-Orlicz norm defined by formula \rf{DO}. In more
detail, for every mapping $T(x)=(P_x,z_x),~ x\in S,$
\be
\|T\|_{LO(S)}&:=&\inf\{\lambda:~\dm(\lambda^{-1}\circ
T(x), \lambda^{-1}\circ T(y))\le \omega(\|x-y\|)~~ {\rm
for~ all}~~ x,y\in S\}\nn\\
&:=& \inf\{\lambda:~\dm((\lambda^{-1}P_x,z_x),
(\lambda^{-1}P_y,z_y))\le \omega(\|x-y\|)~~ {\rm for~
all}~~ x,y\in S\}.
\nn
\ee
\par In Section 1 we have also defined the space
$\BLip(\SO,\TK)$ of all Lipschitz mappings
$T(x)=(P_x,z_x),~x\in S,$ from $\Lip(\SO,\TK)$ such that $
\sup_{x\in S}|D^\alpha P_x(x)|<\infty$ for every $\alpha,
|\alpha|\le k.$ This space is equipped with the ``norm"
$\|\cdot\|^*_{LO(S)}$ defined by \rf{TS}.\medskip
\par {\it Proof of Proposition \ref{WLip}.}
{\it (Necessity.)} Let $F\in\CKO$. We
have to prove that the mapping $T=T(x)=(P_x,x)$, $x\in S,$
where $P_x:=T^k_x(F)$, belongs to $\BLip(\SO,\TK)$. By the
Whitney-Glaeser extension theorem (necessity part)
inequalities  \rf{nd} and \rf{dp} are satisfied with
$\lambda:=c(k,n)\|F\|_{\CKO}$.
\par Put $T_x:=(P_x,x)$. Then inequality \rf{dp} is
equivalent to
$$
\max\{\vfa(\lambda^{-1}|D^\alpha(P_x-P_y)(x)|),
\vfa(\lambda^{-1}|D^\alpha(P_x-P_y)(y)|\}\le\omega(\|x-y\|)
$$
so that
\be &&\dom(\lambda^{-1}\circ T_x,\lambda^{-1}\circ
T_y)\nn\\&& :=\max_{|\alpha|\le k}\{\omega(\|x-y\|),
\vfa(|D^\alpha(\lambda^{-1}P_x-\lambda^{-1}P_y)(x)|),
\vfa(|D^\alpha(\lambda^{-1}P_x-\lambda^{-1}P_y)(y)|\}\nn\\
&&\le\omega(\|x-y\|). \nn \ee
Since $\dm\le\dom$, see Theorem \ref{DOC}, this implies
\bel{TN}
 \dm(\lambda^{-1}\circ T_x,\lambda^{-1}\circ T_y)
\le\omega(\|x-y\|), ~~~x,y\in S,
\ee
which by definition \rf{DO} is equivalent to the inequality
$\|T\|_{LO(S)}\le\lambda$. From this and \rf{nd} we obtain
that $T\in\BLip(\SO,\TK)$ and
$\|T\|^*_{LO(S)}\le 2\lambda=2 c(k,n)\|F\|_{\CKO}. $
\par {\it (Sufficiency).} Assume that the mapping
$T=T(x)=(P_x,x)$, $x\in S,$ belongs to $\BLip(\SO,\TK)$.
Put $\lambda:=2\|T\|^*_{LO(S)}$. Then by \rf{TS} inequality
\rf{nd} of the Whitney-Glaeser extension theorem is
satisfied. Prove that inequality \rf{dp} is true as well.
\par By \rf{DO}
$\|\lambda^{-1}\circ T\|_{\Lip(\SO,\TK)}\le 1$ so that $T$
satisfies inequality \rf{TN}. By Theorem \ref{DOC}
$$
\dom((e^n\lambda)^{-1}\circ T_x,(e^n\lambda)^{-1}\circ T_y)
\le\dm(\lambda^{-1}\circ T_x,\lambda^{-1}\circ T_y)
$$
so that
$$
\dom((e^n\lambda)^{-1}\circ T_x,(e^n\lambda)^{-1}\circ T_y)
\le\omega(\|x-y\|), ~~~x,y\in S.
$$
This inequality and definition \rf{deffi} of $\dom$ imply
that for every $\alpha, |\alpha|\le k,$ and every $x,y\in
S$ we have
\bel{i1} \vfa(|D^\alpha((e^n\lambda)^{-1}P_x-
(e^n\lambda)^{-1}P_y)(x)|)\le\omega(\|x-y\|) \ee
and
\bel{i2} \vfa(|D^\alpha((e^n\lambda)^{-1}P_x-
(e^n\lambda)^{-1}P_y)(y)| \le\omega(\|x-y\|). \ee
Recall that $\vfa:=\omega((s^{k-|\alpha|}\omega(s))^{-1})$
and by our assumption $\omega$ is a strictly increasing
function. This shows that \rf{i1} and \rf{i2} are
equivalent to the required inequality \rf{dp} (with
$e^n\lambda$ instead of $\lambda$.)
\par Thus conditions \rf{nd} and \rf{dp} of the
Whitney-Glaeser extension theorem are satisfied which
implies the existence of a function $F\in\CKO$ such that
$T^k_x(F)=P_x$, $x\in S,$ and $\|F\|_{\CKO}\le
c(k,n)e^n\lambda$. The proposition is proved.\bx
\begin{remark} {\em  Proposition \ref{WLip} allows us to
reformulate the Whitney-Glaeser extension theorem \ref{WG}
as an extension theorem for Lipschitz mappings from
$\BLip(\RN_\omega,\TK)$.
\begin{proposition}
Suppose we are given a family of polynomials $P_x\in\PK,
x\in S,$ such that the mapping $T(x):=(P_x,x),x\in S,$
belongs to $\BLip(\SO,\TK)$. Then $T$ can be extended to a
Lipschitz mapping
$\tT(x)=(\tP_x,x)\in\BLip(\RN_\omega,\TK)$ satisfying
$\|\tT\|^*_{LO(\RN)}\le c(k,n)\|T\|^*_{LO(S)}$.
\end{proposition}
\par {\it Proof.} Since $T\in\BLip(\SO,\TK)$ by Proposition
\ref{WLip} there is a function $F\in\CKO$ such that
$T^k_x(F)=P_x, x\in S,$ and $\|F\|_{\CKO}\le
c(k,n)\|T\|^*_{LO(S)}$. Applying again Proposition
\ref{WLip} (necessity) to the function $F$ (on $\RN$) we
conclude that the mapping $\tT(x):=(T^k_x(F),x), x\in\RN,$
provides the required extension of $T$ from $S$ on all of
$\RN$. Its norm in $\BLip(\RN_\omega,\TK)$ satisfies the
inequality
$$
\|\tT\|^*_{LO(\RN)}\le c_1(k,n)\|F\|_{\CKO}\le
c_1(k,n)c(k,n)\|T\|^*_{LO(S)}.
$$
The proposition is proved.\bx
}\end{remark}
\SECT{3. Lipschitz selections of polynomial-set valued
mappings}{3}
\indent
\par In this section we deal with the Lipschitz
selection problem for the pair of metric spaces
$S_\omega:=(S,\omega(\|\cdot\|))$ and $\TK:=(\PRN,\dm)$.
Our goal is to prove Theorem \ref{FW}. A proof of this
result is based on the classical Helly theorem and a
combinatorial lemma on a structure of finite metric graphs.
For its formulation we let $(\Mc,\rho)$ denote a metric
space. Let $\Tr$ be a (graph-theoretic) tree whose set of
vertices coincides with $\Mc$. If vertices $z,z'$ are
joined by an edge, we write $z\leftrightarrow z'$. This
tree generates a new metric
$$
\rt(x,y):=\sum\limits_{i=0}^{n-1}\rho(z_i,z_{i+1})
$$
where $\{z_0,z_1,...,z_n\}$ is the {\it unique path} in
$\Tr$ joining $x$ and $y$, i.e., $z_0=x, z_n=y, z_i\neq
z_j$ for $i\neq j$ and $z_j$ joined to $z_{j+1}$ by an edge
($z_j\leftrightarrow z_{j+1}$).
\par Clearly, $\rho\le\rt.$ As usual, we let $\dtr x$
denote the degree of a vertex $x$ in $\Tr$, i.e., the
number of edges incident to $x$. Given $a\in\R$, we let
$]a[$ denote an integer $m$ such that $m-1<a\le m$.
\begin{proposition}\label{MT}(\cite{S3})
For every finite metric space $\MR$ there is a tree $\Tr$
such that
$$
\rho(x,y)\le\rt(x,y) \le \eta\rho(x,y), ~~~~~x,y\in\Mc,
$$
and
$$
\max\limits_{x\in\Mc}\dtr x\ge~]\log_2(\card\Mc)[~.
$$
Here $\eta=\eta(\card\Mc)$ is a constant depending only on
cardinality of $\Mc$.
\end{proposition}
\par We turn to\medskip
\par {\it Proof of Theorem \ref{FW}.} Recall that
$\Gc(x)=(G(x),x),~x\in S,$ where $G(x)$ is a convex subset
of $\PK$. Observe also that theorem's statement can be
readily reduced to the case $K=1$. To this end it suffices
to consider a set-valued mapping $\twG(x)=(K^{-1}G(x),x)$
and make use of the fact that given a mapping
$g(x)=(P_x,x), x\in S,$ its norm $\|g\|_{LO(S)}\le K$ iff
$\|\tg\|_{LO(S)}\le 1$ where $\tg(x)=(K^{-1}P_x,x),\, x\in
S.$
\par We prove the theorem by induction on $m:=\card S$.
Put
\bel{LG} \ELG:=\min\{\ell+1,\dim\PK\}. \ee
If $m=2^{\,\ELG}$, nothing to prove. Suppose that the
theorem is true for every set $S$ with $\card S\le m$ where
$m\ge 2^{\,\ELG}$ and prove the result for a set $S$
consisting of $m+1$ points.
\par Thus $\card S=m+1$ and we may assume
that the restriction $\Gc|_{S'}$ to every subset $S'\subset
S$ consisting of at most $m$ points has a Lipschitz
selection $g_{S'}\in\Lip(S'_\omega,\TK)$ such that
$\|g_{S'}\|_{LO(S')}\le 1$. We have to prove that the
set-valued mapping $\Gc$ on all of $S$ has a Lipschitz
selection $g\in\Lip(S_\omega,\TK)$ with $\|g\|_{LO(S)}\le
\gamma(k,n,m)$.
\par Let us apply Proposition \ref{MT} to the metric space
$(S,\rho)$ with $\rho(x,y):=\|x-y\|$. By this proposition
there is a tree $\Tr$ with vertices in $S$ and a vertex
$x_0\in S$  such that $\rho\le\eta(m)\rt$ and
$$
\dtr x_0\ge~ ]\log_2(\card S)[\,\,=\,\,]\log_2 (m+1)[\,\,
\ge\,\,]\log_2(2^{\,\ELG}+1)[\,\,=\ELG+1.
$$
We let $I(x_0)$ denote the family of vertices
$\{y_1,y_2,...,y_p\}$ incident to $x_0$. Thus the number of
these vertices
\bel{pl} p=\dtr x_0\ge \ELG+1. \ee
\par For every vertex $y\in I(x_0)$ we define a subtree
$\Tr_y$ of the tree $\Tr$ whose set of vertices $S_y$
consists of all $z\in S$ for which the (unique) path
connecting $z$ and $y$ in $\Tr$ does not contain the vertex
$x_0$. (We supply $\Tr_y$ with the tree structure induced
by $\Tr$.) Clearly, the trees $\Tr_{y}$ and $\Tr_{y'}$ have
no common vertices for different $y,y'\in I(x_0)$.
\par For each vertex $y\in I(x_0)$ (i.e.,
$y\leftrightarrow x_0$) we let $\Or(y)$ denote a family of
polynomials $P\in\PK$ such that the following holds: For
each vertex $z\in S_y$ of $\Tr_y$ there is a polynomial
$P_z\in G(z)$ such that $P_y=P$ and for every $z,z'\in
S_y$, $z\leftrightarrow z'$, we have
\bel{A3} \dom(e^{-n}\circ T_z, e^{-n}\circ T_{z'})\le
\omega(\|z-z'\|). \ee
where $T_z:=(P_z,z),T_{z'}:=(P_{z'},z')$.
\par Since $\card S_y<\card S$, by the assumption the
restriction $\Gc|_{S_y}$ has a Lipschitz selection
$g_{S_y}:S_y\to \PRN$ with $\|g_{S_y}\|_{LO(S_y)}\le 1$. In
other words, for each $z\in S_y$ there is a polynomial
$P_z\in G(z)$ such that
$$
\dm((P_z,z),(P_{z'},z'))\le \omega(\|z-z'\|).
$$
Then
\bel{S3}
P_y\in\Or(y).
\ee
In fact, by Theorem \ref{DOC}
$$
\dom((e^{-n}P_z,z),(e^{-n}P_{z'},z'))\le
\dm((P_z,z),(P_{z'},z'))\le \omega(\|z-z'\|)
$$
proving \rf{S3}.
\par We have also proved that $\Or(y)\ne\emp$ for every
$y\in I(x_0)$. Recall that inequality \rf{A3} is equivalent
to inequalities \rf{dp} of the Whitney-Glaeser extension
theorem. The left-hand sides of these inequalities are
subadditive and  positively homogeneous functions of
polynomials $P_x, P_y$. This and the definition of $\Or(y)$
show that for every $y\in I(x_0)$ the set $\Or(y)$ is {\it
convex}.
\par Given $y\in I(x_0)$ we put
\be
U(y)&:=&\{P\in\PK: {\rm there~is}~~\tP\in\Or(y)
~{\rm such~that}\nn\\
&&\dom((\ts^{-1}P,x_0),(\ts^{-1}\tP,y))\le
\omega(\|x_0-y\|)\}\label{U3} \ee
where $\ts:=3^{k+1}e^{3n}$.
Prove that
\bel{B3}
G(x_0)\bigcap\{\bigcap_{y\in I(x_0)}U(y)\}\ne\emp.
\ee
But before to do this let us show how the proof of the
theorem can be completed.
\par Property \rf{B3} implies the existence of
polynomials $P_{x_0}\in G(x_0),$ $P_y\in\Or(y)\subset
G(y),$ $y\in I(x_0)$, such that
\bel{C3} \dom((\ts^{-1}P_{x_0},x_0),(\ts^{-1}P_y,y))\le
\omega(\|x_0-y\|). \ee
In turn, since $P_y\in\Or(y)$ for $y\in I(x_0)$, by \rf{A3}
there exist polynomials $P_z\in G(z), z\in S_y$, such that
\bel{D3} \dom((e^{-n}P_z,z),(e^{-n}P_{z'},z'))\le
\omega(\|z-z'\|),~~~z\leftrightarrow z',~ z,z'\in S_y. \ee
\par Now polynomials $P_x$ are defined for all $x\in S$.
Put
$$
g(x):=(P_x,x),~~ x\in S.
$$
Then $g:S\to\TK$ is a selection of $\Gc$. Let us show that
$g\in\Lip(S_\omega,\TK)$ and $\|g\|_{LO(S)}$ is bounded by
a constant depending only on $k,n$ and $m$. In fact, by
\rf{C3} and \rf{D3} for every two vertices $z$ and $z'$ of
the tree $\Tr$ joined by an edge ($z\leftrightarrow z'$) we
have
$$
\dom(\ts^{-1}\circ g(z), \ts^{-1}\circ g(z'))\le
\omega(\|z-z'\|).
$$
By Theorem \ref{DOC} $\dm\le\dom$ so that
\bel{F3}
\dm(\ts^{-1}\circ g(z), \ts^{-1}\circ g(z'))\le
\omega(\|z-z'\|),~~~z\leftrightarrow z', ~z,z'\in S.
\ee
\par To estimate $\dm(g(x),g(y))$ for {\it arbitrary}
$x,y\in S$ we will make use of Lemma \ref{seqest}. Since
$x,y$ are vertices of the tree $\Tr$, there is the unique
path $\{z_0, z_1,...,z_q\}$ in $\Tr$ joining $x$ and $y$
(i.e., $z_0=x$, $z_q=y$ and $z_i\leftrightarrow z_{i+1},$
$i=0,1,...q-1).$ Clearly, $q\le m$ (recall that $\card
S=m+1$).
\par We put $T_i:=\ts^{-1}\circ g(z_i)=
(\ts^{-1}P_{z_i},z_i)$ so that by \rf{F3}
$$
\dm(T_i,T_{i+1})\le\omega(\|z_i-z_{i+1}\|),~~i=0,1,...q-1.
$$
Recall that $\rho(x,y)(:=\|x-y\|)\le \eta\rt(x,y)$ where
$\eta=\eta(m)$ is the constant from Proposition \ref{MT}.
Hence
$$
\sum_{i=0}^{q-1}\|z_i-z_{i+1}\|\le \eta\|z_0-z_q\|~
 (=\eta\|x-y\|).
$$
On the other hand,
\be \sum_{i=0}^{q-1}\omega(\|z_i-z_{i+1}\|)&\le&
q\max_{i=0,...,q-1}\omega(\|z_i-z_{i+1}\|)\le
q\omega\left(\sum_{i=0}^{q-1}\|z_i-z_{i+1}\|\right) \nn\\
&\le& m\,\omega(\eta\|z_0-z_q\|)\le m
\eta\,\omega(\|z_0-z_q\|) \nn \ee
(recall that $\omega$ is a concave non-negative function on
$\R_+$ so that $\omega(\lambda t)\le\lambda \omega(t),$
$\lambda \ge 1$).
\par Let us apply Lemma \ref{seqest} to the family
$\{T_i,~i=0,...,q-1\}$, points $\{z_i,~i=0,...,q-1\}$ and a
parameter $\lambda:=m\eta$. By this lemma
$$
\dm(\tau^{-1}\circ T_0,\tau^{-1}\circ T_q)
\le\omega(\|z_0-z_q\|)
$$
where $\tau:=e^{2n}\lambda^{k+1}$, see \rf{cl}. Since
$z_0=x,z_q=y$ and $T_0=(P_{z_0},z_0)=(P_x,x)=g(x)$,
$T_q=(P_{z_q},z_q)=(P_y,y)=g(y)$, we obtain
$$
\dm(\tau^{-1}\circ g(x),\tau^{-1}\circ g(y))
\le\omega(\|x-y\|),~~~x,y\in S.
$$
Hence $\|g\|_{LO(S)}\le \tau=\tau(k,n,m)$ proving that $g$
is a {\it Lipschitz} selection of $\Gc$.
\par Thus it remains to prove \rf{B3}.
This property readily follows from Helly's theorem and the
induction assumption. We put
\bel{3A}
F(x_0):=G(x_0),~~~~~F(y):=U(y),~~~y\in I(x_0),
\ee
and $\tilde{I}:=\{x_0,y_1,y_2,...,y_p\}\,(={x_0}\cup
I(x_0))$. Then property \rf{B3} is equivalent to
$$
\cap\{F(y):~y\in\tilde{I}\}\ne\emptyset.
$$
By \rf{pl}
$$
\card \tilde{I}=1+\card I(x_0)\ge p+1\ge \ELG+2.
$$
Moreover, all the sets $F(y),y\in\tilde{I},$ are convex
subsets of the finite-dimensional space $\PK$, and
dimension of one of them, of the set $F(x_0):=G(x_0)$, is
at most $\ell$. Therefore by Helly's theorem it suffices to
prove that
$$
\cap\{F(y):~y\in I'\}\ne\emptyset
$$
for every subfamily $I'\subset\tilde{I}$ consisting of at
most
$$
\min\{\ell+2,\dim\PK+1\}=\ELG+1
$$
elements. (Recall that $\ELG$ is defined by \rf{LG}).
\par Since $\card\tilde{I}\ge \ELG+2$ and
$\card I'\le \ELG+1$, there is a point $\ty\in\tilde{I}$
such that $\ty\notin I'$. Then by the assumption for the
set $S':=S\setminus\{\ty\}$ the restriction $\Gc|_{S'}$ has
a Lipschitz selection $g_{S'}:S'\to\TK$ with
$\|g_{S'}\|_{LO(S')}\le 1$. Thus $g_{S'}(y)=(P_y,y), y\in
S',$ where
\bel{PG} P_y\in G(y),~~~ y\in S', \ee
and
\bel{DW} \dm(g_{S'}(y),g_{S'}(y'))
\le\omega(\|y-y'\|),~~~y,y'\in S'. \ee
\par We let $\bar{y}$ denote the nearest to $x_0$
(in the metric $\|\cdot\|$) point from the family $I'$.
(Clearly, $\bar{y}=x_0$ whenever $x_0\in I'$.) Prove that
\bel{3B}
P_{\bar{y}}\in F(y)~~~{\rm for~every}~~y\in I'.
\ee
In fact, if $y=x_0$, then $x_0\in I'$ so that
$\bar{y}=x_0$. Therefore by \rf{3A} $F(y)=G(x_0)$ so that
\rf{3B} follows from \rf{PG}. Thus later on we may assume
that $y\ne x_0$.
\par As we have proved above, see \rf{S3},
$P_y\in\Or(y),y\in I'$. Moreover, by \rf{DW}
$$
\dm(g_{S'}(y),g_{S'}(\bar{y}))
\le\omega(\|y-\bar{y}\|).
$$
On the other hand, by \rf{pdm}
$$
\dm((P_{\bar{y}},\bar{y}),(P_{\bar{y}},x_0))=
\omega(\|\bar{y}-x_0\|).
$$
But by definition of $\bar{y}$
$$
\|x_0-\bar{y}\|+\|\bar{y}-y\|\le 2\|x_0-\bar{y}\|+\|x_0-y\|
\le 3\|x_0-y\|
$$
and
$$
\omega(\|x_0-\bar{y}\|)+\omega(\|\bar{y}-y\|)\le
2\omega(\|x_0-\bar{y}\|)+\omega(\|x_0-y\|) \le
3\omega(\|x_0-y\|).
$$
Now let us apply Lemma \ref{seqest} to
$T_0:=g_{S'}(y)=(P_y,y)$,
$T_1:=g_{S'}(\bar{y})=(P_{\bar{y}},\bar{y})$ and
$T_2:=(P_{\bar{y}},x_0)$ with $\lambda=3$. Then by the
lemma
$$
\dm((3^{k+1}e^{2n})^{-1}\circ T_0,(3^{k+1}e^{2n})^{-1}\circ
T_{2})\le\omega(\|x_0-y\|)
$$
so that by Theorem \ref{DOC}
$$
\dom((3^{k+1}e^{3n})^{-1}\circ
T_0,(3^{k+1}e^{3n})^{-1}\circ
T_{2})\le\dm((3^{k+1}e^{2n})^{-1}\circ
T_0,(3^{k+1}e^{2n})^{-1}\circ T_{2})\le\omega(\|x_0-y\|).
$$
Recall that $\ts:=3^{k+1}e^{3n}$. Hence
$$
\dom(\ts^{-1}\circ T_0,\ts^{-1}\circ
T_{2})=\dom((\ts^{-1}P_y,y),(\ts^{-1}P_{\bar{y}},x_0))
\le\omega(\|x_0-y\|).
$$
But $P_y\in\Or(y)$ so that by definition \rf{U3}
$P_{\bar{y}}\in U(y)=F(y)$ (recall that $y\ne x_0$).
\par Theorem \ref{FW} is completely proved.\bx
\par This theorem implies a similar result for the space
$\BLip(S_\omega,\TK)$ of ``bounded" Lipschitz mappings.
\begin{theorem}\label{BFW} Let $\Gc(x)=(G(x),x),x\in S$,
be a set-valued mapping from a finite set $S\subset\RN$
into $2^{\PRN}$ such that for each $x\in S$ the set
$G(x)\subset\PK$ is a convex set of polynomials of
dimension at most $\ell$. Suppose that for every subset
$S'\subset S$ consisting of at most
$2^{\min\{\ell+1,\dim\PK\}}$ points the restriction
$\Gc|_{S'}$ has a Lipschitz selection
$g_{S'}\in\BLip(S',\TK)$ with $\|g_{S'}\|^*_{LO(S')}\le K$.
Then $\Gc$ on all of $S$ has a Lipschitz selection
$g\in\BLip(S,\TK)$ with
 $\|g\|^*_{LO(S)}\le \gamma(k,n,\card S)K$.
\end{theorem}
\par {\it Proof.} As in the proof of Theorem \ref{FW} it
suffices to prove the result for $K=1$. Given $x\in S$ we
put
$$
H(x):=\{P\in\PK:~\max_{|\alpha|\le k}|D^\alpha P_x(x)| \le
1\}.
$$
We define a set-valued mapping $\twG$ by letting
\bel{B4}
\twG(x):=(G(x)\cap H(x),x), ~~~x\in S.
\ee
\par Put $\ELG:=\min\{\ell+1,\dim\PK\}$ and prove that for
every subset $S'\subset S$ of cardina\-lity $\card S'\le
2^{\ELG}$ the restriction $\twG|_{S'}$ has a Lipschitz
selection $\tg_{S'}\in\Lip(S'_\omega,\TK)$ with
$\|\tg_{S'}\|_{LO(S')}\le 1$. In fact, by theorem's
hypothesis $\Gc|_{S'}$ has a selection
$g_{S'}\in\Lip(S'_\omega,\TK)$ such that
$\|g_{S'}\|^*_{LO(S')}\le 1$. Thus
$g_{S'}(x)=(P_{(S',x)},x)$ where the polynomial
$P_{(S',x)},\, x\in S',$ satisfy the following conditions:
(i). $P_{(S',x)}\in G(x), x\in S';$~ (ii). $|D^\alpha
P_{(S',x)}(x)|\le 1$ for all $|\alpha|\le k$ and $x\in S'$,
and (iii).
\bel{A4}
\dm((P_{(S',x)},x),(P_{(S',y)},y))\le\omega(\|x-y\|),
~~~x,y\in S'. \ee
Hence $P_{(S',x)}\in H(x)\cap G(x), x\in S',$ so that the
mapping $\tg_{S'}(x):=(P_{(S',x)},x), x\in S',$ provides
the required selection of $\twG|_{S'}$. By \rf{A4} its
Lipschitz-Orlicz norm in  $\Lip(S'_\omega,\TK)$ does not
exceed $1$.
\par By Theorem \ref{FW} $\twG$ on all of $S$ has a
Lipschitz selection $g(x):=(P_x,x), x\in S,$ with
$\|g\|_{LO(S)}\le \gamma_1(k,n,\card S)$. Since $g$ is a
selection of $\twG$, by \rf{B4} it is a selection of $\Gc$
as well. Moreover, by \rf{B4} $P_x\in H(x), x\in S,$ so
that $\max_{|\alpha|\le k}|D^\alpha P_x(x)|\le 1$ for all
$x\in S$. Hence
$$
\|g\|^*_{LO(S')}=\max_{|\alpha|\le k}\, \sup_{x\in
S}|D^\alpha P_x(x)|+\|g\|_{LO(S)}\le 1+\gamma_1(k,n,\card
S).
$$
\par The theorem is proved.\bx
\SECT{4. The weak finiteness property of the space
$\CKO$}{4}
\indent
\par {\it Proof of Theorem \ref{TF1}.} The result easily
follows from Proposition \ref{WLip} and Theorem \ref{BFW}.
In fact, we let $\Gc$ denote a set-valued mapping
$\Gc(x):=(G(x),x), x\in S.$
\par Fix a set $S'\subset S$ of cardinality at most
$2^{\,\ELG}$ where $\ELG:=\min\{\ell+1,\dim\PK\}$. By
theorem's hypothesis there is a function $F_{S'}\in \CKO$
with $\|F_{S'}\|_{\CKO}\leq 1$ satisfying
\bel{E4} T_{x}^{k}(F_{S'})\in G(x), ~~~~ x\in S'. \ee
We put $T_{S'}(x):=(T_{x}^{k}(F_{S'}),x),\,x\in S'.$ By
\rf{E4} $T_{S'}$ is a selection of the restriction
$\Gc|_{S'}$. Moreover, by Proposition \ref{WLip} (``only
if" part) the mapping $T_{S'}:S'\to \TK$ belongs to
$\BLip(S_\omega,\TK)$ and its Lipschitz-Orlicz norm
satisfies the inequality
$$
\|T_{S'}\|^*_{LO(S')}\le c_1(k,n)\|F_{S'}\|_{\CKO}\le
c_1(k,n).
$$
Thus $T_{S'}$ is a Lipschitz selection of $\Gc|_{S'}$.
Since $S'$ is an arbitrary subset of $S$ of cardinality at
most $2^{\,l_G}$, by Theorem \ref{BFW} there is a selection
$T(x)=(P_x,x)$ of $\Gc$ defined on all of $S$ and
satisfying $\|T\|^*_{LO(S')}\le c_1(k,n)\gamma(k,n,\card
S)$.
\par In particular, $P_x\in G(x),\, x\in S$. Now by
Proposition \ref{WLip} (``if" part) there is a function
$F\in\CKO$ with
$$
\|F\|_{\CKO}\le c_2(k,n)\|T\|^*_{LO(S')}\le
c_2(k,n)c_1(k,n)\gamma(k,n,\card S)
$$
such that $T_{x}^{k}(F)=P_x,\, x\in S.$ Hence
$T_{x}^{k}(F)\in G(x), x\in S,$ and the theorem follows.\bx
\begin{remark}{\em Theorem \ref{TF3} for
$\xi\equiv 0$ implies the following finiteness property of
the space $\CKO$: {\it A function $f$ defined on a subset
$S\subset\RN$ can be extended to a function $F\in \CKO$
with $\|F\|_{\CKO}\le \gamma(k,n)$ provided its restriction
$f|_{S'}$ to every subset $S'\subset S$ consisting of at
most $N(k,n)=2^{\dim\PK}$ points can be extended to a
function $F_{S'}\in \CKO$ with $\|F_{S'}\|_{\CKO}\le 1$.}
\par In particular, $\dim \Pc_1=n+1$ so
that $N(1,n)\le 2^{n+1}$. Recall that the sharp value of
the finiteness number for $k=1$ equals $3\cdot
2^{n-1}=\frac{3}{4}\cdot 2^{n+1}$. This shows that the
estimate $2^{\dim\PK}$ is rather far from the optimal one
and apparently can be decreased considerably.
\par In the next paper we will prove that the finiteness
number $N(k,n)$ does not exceed $(k+1)\cdot 2^{\dim\PK-k}$.
We conjecture that the sharp value of the finiteness number
in the above finiteness property for $\CKO$ is
$$
N(k,n)=\prod_{m=0}^{k}(k-m+2)^{{n+m-2\choose m}}.
$$
}\end{remark}
\renewcommand {\refname} {\centerline{\normalsize
{\bf References}}}

\end{document}